\documentclass[12pt,a4paper]{article}
\usepackage{amssymb,amsmath,amsthm,euscript}
\usepackage{leftidx}
\usepackage[matrix,arrow,curve]{xy}

\usepackage[utf8]{inputenc}
\usepackage[english]{babel}

\newcommand{\chara}{\mathop{\mathrm{char}}\nolimits}
\newcommand{\End}{\mathop{\mathrm{End}}\nolimits}

\newcommand{\Hom}{\mathop{\mathrm{Hom}}\nolimits}

\newcommand{\tr}{\mathop{\mathrm{tr}}\nolimits}

\newcommand{\rk}{\mathop{\mathrm{rk}}\nolimits}

\newcommand{\Min}{\mathop{\mathrm{Min}}\nolimits}

\numberwithin{equation}{section}

\renewcommand{\le}{\leqslant}
\renewcommand{\ge}{\geqslant}

\newcommand{\la}{\langle}
\newcommand{\ra}{\rangle}

\newcommand{\B}{{\mathcal B}}

\newcommand{\gR}{{\mathfrak R}}

\newcommand{\gX}{{\mathfrak X}}
\newcommand{\gU}{{\mathfrak U}}

\newcommand{\bfq}{{\mathbf q}}

\newcommand{\wt}{\widetilde}
\newcommand{\wh}{\widehat}

\newtheorem{Th}{Theorem} [section]
\newtheorem{Lem}[Th]{Lemma}

\newcommand{\KK}{{\mathbb{K}}}

\newcommand{\SSS}{{\mathbb{S}}}

\usepackage{hyperref}
\newcommand{\eprint}[1]{\href{https://arxiv.org/abs/#1}{\tt arXiv:#1}}

\newcounter{bbcount}[subsection]

\title{Manin matrices of type $C$: multi-parametric deformation}
\author{Alexey Silantyev\thanks{aleksejsilantjev@gmail.com}}
\date{}

\begin{document}

\maketitle

\vspace{-5mm}
\begin{center}
{\it Joint Institute for Nuclear Research, 141980 Dubna, Moscow region, Russia} \\
{\it State University ``Dubna''{}, 141980 Dubna, Moscow region, Russia} \\
\end{center}
\vspace{5mm}

\begin{abstract}
 We constructed a multi-parametric deformation of the Brauer algebra representation related with the symplectic Lie algebras. The notion of Manin matrix of type $C$ was generalised to the case of the multi-parametric deformation by using this representation and corresponding quadratic algebras. We derived pairing operators for these quadratic algebras and minors for the considered Manin matrices. The rank of pairing operators and dimensions of components of quadratic algebras were calculated. 
\end{abstract}


\tableofcontents

\section{Introduction}

In~\cite{Manin87,Manin88} Yuri Manin considered some `non-commutative' transformations between arbitrary quadratic algebras. He interpreted the quadratic algebras and these transformations as a quantum version of vector spaces and representations on them (see~\cite{Sqrt} for details).

The notion `Manin matrix' was introduced in~\cite{CF} for a matrix over a non-commutative ring satisfying some commutation relations. Such a matrix describes a non-commutative transformation between polynomial algebras. The Manin matrices and their determinants have a lot of interesting properties and applications~\cite{CF,CM,CFR}. In~\cite{CFRS} they were generalised for the $q$-case, where $q$ is a parameter of the deformation of the polynomial algebras. The MacMahon Master Theorem and Cayley--Hamilton--Newton identities for the $q$-Manin matrices were obtained in~\cite{GLZ,FH07,FH08} and \cite{IO} respectively. Some properties and applications of the super-version of Manin matrices were derived in~\cite{MR}.

The Manin matrices were described for arbitrary quadratic algebras in terms of idempotent operators in~\cite{Smm}. In this general case we cannot generalise the notion of determinant, but we introduced two types of minors for Manin matrices by using non-degenerate pairings between quadratic algebras.

Multi-parametric deformation of the polynomial algebras and corresponding matrices were considered by Manin for the super-case in~\cite{Manin89}. In~\cite[\S~3.3, \S~6.1]{Smm} the (non-super) multi-parametric case was considered: we obtained idempotents describing the corresponding quadratic algebras and calculated the minors of Manin matrices by using representations of the symmetric groups. In terms of the classification of simple Lie algebras one can refer to this case as `multi-parametric Manin matrices of type $A$'.

Manin matrices of types $B$, $C$ and $D$ were introduced by A. Molev in~\cite{MolevSO} for the non-deformed case. They are related with the representations of Brauer algebras. The corresponding quadratic algebras and minors of Manin matrices of these types were described in \cite[\S~7]{Smm}. Since the quadratic algebras over an algebraically closed field of characteristic $0$ for different orthogonal/symplectic forms are isomorphic to each other, it is enough to take a canonical form in this case.

Here we generalise the Manin matrices of type $C$ for the multi-parametric case. We do not suppose that the basic field is algebraically closed, but any symplectic form is reduced to the canonical form over an arbitrary field of characteristic $0$. For the orthogonal case (types $B$ and $D$) one needs to consider a general symmetric non-degenerate form, we leave this case for further research.

The structure of the paper is following. In Section~\ref{secMM} we recall necessary notions and facts on quadratic algebras and Manin matrices. In Section~\ref{secA} we recall the multi-parametric type $A$ case. We also recall some definitions and statements on Brauer algebra in Section~\ref{secId}. In Section~\ref{secRep} we construct a deformed representation of the Brauer algebra.  In Section~\ref{secC} we use this representation to define a multi-parametric version of Manin matrices of type $C$ and write some of their minors. In Section~\ref{secDim} we calculate dimensions of the graded components of the quadratic algebras (which equal to the ranks of pairing operators) and prove that the other minors vanish.

\vspace{3mm}
{\it Acknowledgements}. 
The author thanks V. Rubtsov, A. Isaev and anonymous referee for useful references and advice.

\section{Manin matrices for quadratic algebras}
\label{secMM}

We fix a field $\KK$ of characteristic $\chara\KK=0$ and an integer $k\ge2$. Recall the general notions introduced in~\cite{Smm}.

Let $V$ be a finite-dimensional vector space (over $\KK$) with a basis $(e_i)$ labelled by some set of indices. Denote by $E_i^j\in\End(V)$ the operators acting by the formula $E_i^je_l=\delta^j_{l}e_i$. They form a basis of $\End(V)$. We use the following tensor notations. Let $T$ be an operator acting on the vector space $V\otimes V=V\otimes_\KK V$. It has the form $T=\sum_{i,j,l,m}T^{ij}_{lm}E_i^l\otimes E_j^m$, where $T^{ij}_{lm}\in\KK$ are its entries in the basis $(e_i\otimes e_j)$. For $a,b\in\{1,\ldots,k\}$ such that $a\ne b$ denote by $T^{(ab)}$ or $T^{(a,b)}$ the following operator acting on the vector space $V^{\otimes k}=V\otimes\cdots\otimes V$:
\begin{align}
 T^{(ab)}=\sum_{i,j,l,m}T^{ij}_{lm}(1\otimes\cdots\otimes E_i^l\otimes\cdots\otimes E_j^m\otimes\cdots\otimes1),
\end{align}
where $E_i^l$ and $E_j^m$ are on the $a$-th and $b$-th sites respectively.

Let $\wt V$ be a finite-dimensional vector space with a basis $(e_\alpha)$. The vector space $\Hom(\wt V,V)$ has a basis consisting of the operators $E^\alpha_i$ acting as $E^\alpha_i e_\beta=\delta^\alpha_\beta e_i$.
Let $\gR$ be an algebra (i.e. associative unital algebra over $\KK$) and let $M=\sum_{i,\alpha}M_\alpha^i\otimes E_i^\alpha\in\gR\otimes\End(V)$. We identify the elements $r\in\gR$ and $H\in\Hom(\wt V,V)$ with the elements $r\otimes1$ and $1\otimes H$ of $\gR\otimes\Hom(\wt V,V)$ respectively, so that $M=\sum_{i,\alpha}M_\alpha^i E_i^\alpha$. Denote by $M^{(a)}$ an element of $\gR\otimes\Hom(\wt V^{\otimes k},V^{\otimes k})=\gR\otimes\Hom(\wt V,V)\otimes\cdots\otimes\Hom(\wt V,V)$ of the form
\begin{align}
M^{(a)}=\sum_{i,\alpha}M_\alpha^i(1\otimes\cdots\otimes E_i^\alpha\otimes\cdots\otimes1),
\end{align}
where $E_i^\alpha$ is in the $a$-th tensor factor $\Hom(\wt V,V)$.

Let $(e^i)$ be a basis of the vector space $V^*$ dual to the basis $(e_i)$. Let $A\in\End(V\otimes V)$ be an idempotent operator, i.e. $A^2=A$. It acts on the basis vectors $e_i\otimes e_j\in V\otimes V$ and $e^i\otimes e^j\in V^*\otimes V^*$ from the left and from the right respectively: $A(e_i\otimes e_j)=\sum_{l,m}A^{lm}_{ij}e_l\otimes e_m$, $(e^i\otimes e^j)A=\sum_{l,m}A_{lm}^{ij}e^l\otimes e^m$. Denote by $\gX_A(\KK)$, $\gX^*_A(\KK)$, $\Xi_A(\KK)$ and $\Xi^*_A(\KK)$ the quadratic algebras generated by $x^i$, $x_i$, $\psi_i$, $\psi^i$ respectively with the relations
\begin{align} \label{AXX}
 &A(X\otimes X)=0, &&(X^*\otimes X^*)A=0,
 &&(\Psi\otimes\Psi)(1-A)=0, &&(1-A)(\Psi^*\otimes\Psi^*)=0,
\end{align}
where $X=\sum_i x^ie_i$, $X^*=\sum_i x_ie^i$, $\Psi=\sum_i \psi_i e^i$, $\Psi^*=\sum_i \psi^ie_i$.

In~\cite[\S~5.3]{Smm} we introduce the notion of pairing operators $S_{(k)},A_{(k)}\in\End(V^{\otimes k})$ for an idempotent $A$. These operators define non-degenerate pairings $\gX_A(\KK)\times\gX^*_A(\KK)\to\KK$ and $\Xi^*_A(\KK)\times\Xi_A(\KK)\to\KK$ by the matrix formulae
\begin{align}
 &\la X\otimes\cdots\otimes X,X^*\otimes\cdots\otimes X^*\ra=S_{(k)},
 &&\la \Psi^*\otimes\cdots\otimes\Psi^*,\Psi\otimes\cdots\otimes\Psi\ra=A_{(k)},
\end{align}
where the dots `$\cdots$' mean that there are $k$ tensor factors in the products. Recall also that $A_{(1)}=S_{(1)}=1$, $A_{(2)}=A$, $S_{(2)}=1-A$.

Consider two idempotents $A\in\End(V\otimes V)$ and $\wt A\in\End(\wt V\otimes\wt V)$ with entries $A^{ij}_{lm}$, $\wt A^{\alpha\beta}_{\gamma\delta}$. Recall that {\it Manin matrix} over an algebra $\gR$ for a pair of idempotents $(A,\wt A)$ or simply {\it $(A,\wt A)$-Manin matrix} is an operator $M\in\gR\otimes\Hom(\wt V,V)$ with entries $M^i_\alpha\in\gR$ satisfying the relation
\begin{align}
 AM^{(1)}M^{(2)}(1-\wt A)=0.
\end{align}
Entry-wise it is written as $\sum_{l,m}A^{ij}_{lm}M^l_\gamma M^m_\delta=\sum_{l,m,\alpha,\beta}A^{ij}_{lm}M^l_\alpha M^m_\beta\wt A^{\alpha\beta}_{\gamma\delta}$. The Manin matrices are related with the quadratic algebras as follows. There is a bijection between $(A,\wt A)$-Manin matrices $M$ over $\gR$, graded homomorphisms $f_M\colon\gX_{A}(\KK)\to\gR\otimes\gX_{\wt A}(\KK)$ and graded homomorphisms $f^M\colon\Xi_{\wt A}(\KK)\to\gR\otimes\Xi_A(\KK)$ (see~\cite[\S~2.5]{Smm} for details).

In~\cite[\S~5.4]{Smm} we defined two generalisations of minors: $S$-minors and $A$-minors. For an $(A,\wt A)$-Manin matrix $M$ they are entries of the operators
\begin{align}
 &\Min_{\wt S_{(k)}}M=M^{(1)}\cdots M^{(k)}\wt S_{(k)}=\la f_M(X\otimes\cdots\otimes X),\wt X^*\otimes\cdots\otimes\wt X^*\ra, \label{MinS} \\
 &\Min^{A_{(k)}}M=A_{(k)}M^{(1)}\cdots M^{(k)}=\la\Psi^*\otimes\cdots\otimes\Psi^*,f^M(\wt\Psi\otimes\cdots\otimes\wt\Psi)\ra, \label{MinA}
\end{align}
where $\wt S_{(k)}$ and $A_{(k)}$ are pairing $S$- and $A$-operators for the idempotents $\wt A$ and $A$ respectively (see~\cite[\S~5.3]{Smm}). The $S$- and $A$-minors of $M$ describe the $k$-th graded component of the homomorphism $f_M$ and $f^M$ respectively.

\section{Manin matrices of type $A$}
\label{secA}

Here we recall multi-parametric case described in~\cite{Manin89} and~\cite[\S~3.3, 6.1]{Smm}.

Recall that the permutation operator $P=\sum_{i,j}E_i^j\otimes E_j^i$ acting on the space $V\otimes V$ defines two representations of the permutation group $\SSS_k$ on the space $V^{\otimes k}$, these are $\sigma_a\mapsto\pm P^{(a,a+1)}$, where $\sigma_a=\sigma_{a,a+1}\in\SSS_k$ are adjacent transpositions. The polynomial algebra and Grassmann algebra are quadratic algebras $\gX_A(\KK)$ and $\Xi_A(\KK)$ for the idempotent $A=\frac{1-P}2$.

Let $\wh q$ be a matrix with non-zero entries $q_{ij}\in\KK\backslash\{0\}$. We call it a {\it parameter matrix} iff
\begin{align} \label{qpm}
 &q_{ij}=q_{ji}^{-1}, &&q_{ii}=1.
\end{align}
Define ${\wh q}$-permutation operator
\begin{align} \label{Pq}
 P_{\wh q}=\sum_{i,j}q_{ij}E_j^i\otimes E_i^j.
\end{align}
This is an operator $P_{\wh q}\in\End(V\otimes V)$ with the entries $(P_{\wh q})^{lm}_{ij}=q_{ij}\delta^l_j\delta^m_i$.
It is involutive and satisfies the braid relation:
\begin{align} \label{Braid}
 &P_{\wh q}^2=1, &&P_{\wh q}^{(12)}P_{\wh q}^{(23)}P_{\wh q}^{(12)}=P_{\wh q}^{(23)}P_{\wh q}^{(12)}P_{\wh q}^{(23)}.
\end{align}
Due to~\eqref{Braid} the formulae $\sigma_a\mapsto\pm P_{\wh q}^{(a,a+1)}$ gives representations $\rho^\pm_{\wh q}\colon\SSS_k\to\End(V^{\otimes k})$ of the permutation group $\SSS_k$.

The idempotent $A_{\wh q}=\frac{1-P_{\wh q}}2$ defines multi-parametric deformations of the polynomial and Grassmann algebras. These are the quadratic algebras $\gX_{A_{\wh q}}(\KK)$ and $\Xi_{A_{\wh q}}(\KK)$ defined by the commutation relations $x^ix^j=q_{ij}^{-1}x^jx^i$ and $\psi_i\psi_j=-q_{ij}\psi_j\psi_i$. Since the relations for $\gX^*_{A_{\wh q}}(\KK)$ and $\Xi^*_{A_{\wh q}}(\KK)$ are $x^ix^j=q_{ij}x^jx^i$ and $\psi^i\psi^j=-q^{-1}_{ij}\psi^j\psi^i$, we have isomorphisms $\gX^*_{A_{\wh q}}(\KK)\cong\gX^*_{A_{\wh q'}}(\KK)$ and $\Xi^*_{A_{\wh q}}(\KK)\cong\Xi^*_{A_{\wh q'}}(\KK)$, where ${\wh q}\,'=(q_{ij}^{-1})$. The pairing $S$- and $A$-operators for $A_{\wh q}$ are
\begin{align}
 &S_{\wh q,(k)}=\rho^+_{\wh q}(h_{(k)}), &
 &A_{\wh q,(k)}=\rho^-_{\wh q}(h_{(k)}),
\end{align}
where $h_{(k)}=\frac1{k!}\sum_{\sigma\in\SSS_k}\sigma$ and $\rho^\pm_{\wh q}\colon\KK[\SSS_k]\to\End(V^{\otimes k})$ are the corresponding representations of the group algebra $\KK[\SSS_k]$.

Consider an $\wt n\times\wt n$ parameter matrix $\wh p$, where $\wt n=\dim\wt V$. An $(A_{\wh q},A_{\wh p})$-Manin matrix is called briefly {\it $(\wh q,\wh p)$-Manin matrix}. Its $S$- and $A$-minors are defined by the formulae~\eqref{MinS} and \eqref{MinA} with $\wt S_{(k)}=S_{\wh p,(k)}$ and $A_{(k)}=A_{\wh q,(k)}$.

Finally recall some explicit formulae needed below. We have the following relations in the algebras $\Xi_{A_{\wh q}}(\KK)$ and $\Xi^*_{A_{\wh q}}(\KK)$: for any $k$-tuple of indices $I=(i_1,\ldots,i_k)$ and permutation $\sigma\in\SSS_k$ we have
\begin{align}
 \psi_{i_{\sigma(1)}}\cdots\psi_{i_{\sigma(k)}}&=\frac{(-1)^\sigma}{\mu_I(\wh q,\sigma)}\psi_{i_1}\cdots\psi_{i_k}, &
 \psi^{i_{\sigma(1)}}\cdots\psi^{i_{\sigma(k)}}&=(-1)^\sigma\mu_I(\wh q,\sigma)\psi^{i_1}\cdots\psi^{i_k}, \label{Lempsisigma2k}
\end{align}
where $(-1)^\sigma$ is the sign of $\sigma$ and
\begin{align} \label{muI}
&\mu_I(\wh q,\sigma)
=\prod\limits_{s<t\atop\sigma^{-1}(s)>\sigma^{-1}(t)}q_{i_si_t}.
\end{align}
The non-vanishing entries of the $A$-operator $A_{\wh q,(k)}$ are
\begin{align} \label{Awhq}
 (A_{\wh q})^{i_{\sigma(1)}\ldots i_{\sigma(k)}}_{i_{1} \ldots i_{k}}=
 \la\psi^{i_{\sigma(1)}}\cdots\psi^{i_{\sigma(k)}},\psi_{i_{1}}\cdots\psi_{i_{k}}\ra=\frac1{k!}(-1)^\sigma\mu_I(\wh q,\sigma),
\end{align}
where $\sigma\in\SSS_k$ and $I=(i_1,\ldots,i_k)$ is a $k$-tuple of pairwise different indices%
% footnote
\footnote{In~\cite[eq.~(6.5)]{Smm} we proved the formula~\eqref{Awhq} for the case $i_1<\ldots<i_k$, but it still holds for any pairwise different indices since it does not depend on the order in the set of indices. Explicitly one can check~\eqref{Awhq} by using the formula $\mu_I(\wh q,\sigma\tau)=\mu_J(\wh q,\tau)\mu_I(\wh q,\sigma)$, where $J=(i_{\sigma(1)},\ldots,i_{\sigma(k)})$, which in turn follows from~\eqref{Lempsisigma2k}.
}
% end of footnote
.

\section{Idempotents in Brauer algebra}
\label{secId}

We introduce the Brauer algebra~\cite{Br} by following~\cite{Nazarov,BW} (see also~\cite{IM,IMO,MolevSO}). Let $\omega\in\KK$ be a parameter. The Brauer algebra $\B_k(\omega)$ is an algebra generated by the elements $\sigma_1,\ldots,\sigma_{k-1}$ and $\epsilon_1,\ldots,\epsilon_{k-1}$ with the commutation relations
\begin{gather}
 \sigma_a^2=1, \qquad \epsilon_a^2=\omega\epsilon_a, \qquad \sigma_a\epsilon_a=\epsilon_a\sigma_a=\epsilon_a, \qquad\qquad a=1,\ldots,k-1, \label{DefBr1} \\
 \sigma_a\sigma_b=\sigma_b\sigma_a, \qquad \epsilon_a\epsilon_b=\epsilon_b\epsilon_a, \qquad \sigma_a\epsilon_b=\epsilon_b\sigma_a, \qquad\qquad|a-b|>1, \label{DefBr2} \\
 \sigma_a\sigma_{a+1}\sigma_a=\sigma_{a+1}\sigma_a\sigma_{a+1}, \qquad
 \epsilon_a\epsilon_{a+1}\epsilon_a=\epsilon_a, \qquad
 \epsilon_{a+1}\epsilon_a\epsilon_{a+1}=\epsilon_{a+1},\qquad\qquad\qquad\qquad \notag \\
 \sigma_a\epsilon_{a+1}\epsilon_a=\sigma_{a+1}\epsilon_a, \qquad
 \epsilon_{a+1}\epsilon_a\sigma_{a+1}=\epsilon_{a+1}\sigma_a, \qquad\qquad a=1,\ldots,k-2. \label{DefBr3}
\end{gather}
The subalgebra of $\B_k(\omega)$ generated by the elements $\sigma_1,\ldots,\sigma_{k-1}$ is naturally identified with the group algebra $\KK[\SSS_k]$.

Suppose $\omega\ne0$. The relations~\eqref{DefBr1} imply that the elements
\begin{align} \label{va}
 &v_a=\frac{1+\sigma_a}2-\frac{\epsilon_a}{\omega}, &&a=1,\ldots,k-1,
\end{align}
are idempotents: $v_a^2=v_a$.

In order to consider `higher' idempotents we need the elements
 $\sigma_{ab}=\sigma_{ba}=\tau_{a,b-1}\sigma_{b-1}\tau_{a,b-1}^{-1}$, 
 $\epsilon_{ab}=\epsilon_{ba}=\tau_{a,b-1}\epsilon_{b-1}\tau_{a,b-1}^{-1}$, $a<b$,
where $\tau_{ab}=\sigma_a\sigma_{a+1}\cdots\sigma_{b-1}$ is the cycle of the length $b-a$ in the group $\SSS_k$. In particular, $\tau_{aa}=1$, so we have $\sigma_{a,a+1}=\sigma_{a+1,a}=\sigma_a$ and $\epsilon_{a,a+1}=\epsilon_{a+1,a}=\epsilon_a$. Note that the elements $\sigma_{ab}\in\SSS_k$ are the standard transpositions and that $\sigma\epsilon_{ab}\sigma^{-1}=\epsilon_{\sigma(a),\sigma(b)}$ for any $\sigma\in\SSS_k$.

The idempotents were obtained in~\cite{HS,I} for the more general Birman--Murakami--Wenzl algebra. In the particular case of the Brauer algebra the expression found in~\cite{I} reduces as follows (see~\cite{IR,MolevSO,Smm}):
\begin{align}
 &s_{(k)}=\frac1{k!}\prod_{b=2}^k\frac{(y_b+1)(y_b+\omega+b-3)}{2b+\omega-4}, \label{sk}
\end{align}
where $y_b=\sum_{a=1}^{b-1}(\sigma_{ab}-\epsilon_{ab})$ and we suppose that the denominators do not vanish:
\begin{align} \label{omegane}
 \omega\notin\{0,-2,-4,\ldots,4-2k\}. 
\end{align}
For $k=2$ the idempotent~\eqref{sk} coincides with~\eqref{va}, i.e. $s_{(2)}=v_1$. An `$R$-matrix' form of the idempotents $s_{(k)}$ was found in~\cite{IM,IMO}.

It was proved in~\cite[\S~4.3]{I} (see also~\cite{HS}, \cite[eq.~(1.31)]{MolevSO}, \cite[Prop.7.10]{Smm}) that~\eqref{sk} satisfies
\begin{align} \label{sigmask}
 &\sigma_as_{(k)}=s_{(k)}\sigma_a=s_{(k)}, &
 &\epsilon_as_{(k)}=s_{(k)}\epsilon_a=0,  &&a=1,\ldots,k-1.
\end{align}

\section{Multi-parametric representation of Brauer algebra}
\label{secRep}

Now let us consider the case associated with the symplectic Lie algebra $\mathfrak{sp}_{2r}(\KK)$, where $r$ is the rank. In this case we suppose that $V$ has dimension $2r$ and the parameter of the Brauer algebra is $\omega=-2r$. Let the basis $(e_i)$ be labelled by the set $\{-r,\ldots,-1,+1,\ldots,r\}.$
Consider $Q=\sum_{i,j}\varepsilon_i\varepsilon_j E_i^j\otimes E_{-i}^{-j}$, where $\varepsilon_i=1$ for $i>0$ and $\varepsilon_i=-1$ for $i<0$.
The formulae
\begin{align} \label{RepCan}
 &\sigma_a\mapsto -P^{(a,a+1)}, &&\sigma_a\mapsto -Q^{(a,a+1)} 
\end{align}
define a representation of the Brauer algebra $\B_k(-2r)$ on the vector space $V^{\otimes k}$ (see e.g.~\cite{Br,Nazarov,IMO,MolevSO}).
 This representation is related with the canonical anti-diagonal symplectic form $\Omega_{ij}=\varepsilon_i\delta_{i,-j}$. The representation related with the general symplectic form $\Omega=(\Omega_{ij})$ is given by the operator $Q=\sum_{ij}\Omega^{ij}\Omega_{ml}E^l_i\otimes E^m_j$, where $(\Omega^{ij})$ is the matrix inverse to $\Omega$; this representation is equivalent to~\eqref{RepCan}. We construct a multi-parametric deformation of the representation~\eqref{RepCan} by extending the representation $\rho^-_{\wh q}\colon\KK[\SSS_k]\to\End(V^{\otimes k})$ for a parameter matrix $\wh q=(q_{ij})$.

Consider a deformation of $Q$ of the form
\begin{align} \label{Qq}
 &Q_{\vec q}=\sum_{i,j}\varepsilon_i\varepsilon_j q_iq_j^{-1} E_i^j\otimes E_{-i}^{-j}, &
&(Q_{\vec q})^{il}_{jm}=\varepsilon_i\varepsilon_jq_iq_j^{-1}\delta^{i,-l}\delta_{j,-m},
\end{align}
where $\vec q=(q_i)$ is a row-vector with entries $q_i\in\KK\backslash\{0\}$ such that
\begin{align} \label{vecq}
 &q_{-i}=q_i^{-1}, &&q_{-i,j}=q_{ij}^{-1}q_j^2.
\end{align}
In particular, $q_{-i,i}=q_i^2$, so $\vec q$ is defined by $\wh q$ up to signs of $q_i$ (the substitution $q_i\to-q_i$ corresponds to the renumbering of the basis $e_i\leftrightarrow e_{-i}$). Let $\bfq=(\wh q,\vec q)$ be a family of the parameters $q_{ij},q_i$ satisfying~\eqref{qpm}, \eqref{vecq}. It is defined by $r(r+1)/2$ independent parameters: $q_1,\ldots,q_r$, $q_{ij}$, $1\le i<j\le r$. Note that the simultaneous change $q_i\to-q_i$ preserves the operator~\eqref{Qq}. Sometimes we use the notations $P_\bfq=P_{\wh q}$ and $Q_\bfq=Q_{\vec q}$.

\begin{Th} A family $\bfq=(q_{ij},q_i)$ satisfying~\eqref{qpm} and \eqref{vecq} defines a representation $\rho_\bfq\colon\B_k(-2r)\to\End(V^{\otimes k})$ by the formulae
\begin{align} \label{Rep}
 &\sigma_a\mapsto -P_\bfq^{(a,a+1)}, &&\epsilon_a\mapsto -Q_\bfq^{(a,a+1)}, &&a=1,\ldots,k-1.
\end{align}
\end{Th}

\noindent{\bf Proof.} We already know that the commutation relations between $\rho_\bfq(\sigma_a)$ are satisfied. The last two relations~\eqref{DefBr1} follow from
\begin{align}
 &Q_\bfq^2=2r Q_\bfq, &&P_\bfq Q_\bfq = Q_\bfq P_\bfq =-Q_\bfq.
\end{align}
The relations~\eqref{DefBr2} are obvious. Then, it is straightforward to check
\begin{align}
 &Q_\bfq^{(12)}Q_\bfq^{(23)}Q_\bfq^{(12)}=Q_\bfq^{(12)}, &
 &P_\bfq^{(12)}Q_\bfq^{(23)}Q_\bfq^{(12)}=-P_\bfq^{(23)}Q_\bfq^{(12)}. \label{QQPQ}
\end{align}
These imply the second and forth relations~\eqref{DefBr3}. Let $\bfq'$ be the family of the inverted parameters: $\bfq'=(q_{ij}^{-1},q_i^{-1})$. We have
\begin{align} \label{PQ21}
 &P_{\bfq'}=P_\bfq^{(21)}, &&Q_{\bfq'}=Q_\bfq^{(21)}.
\end{align}
By applying~\eqref{PQ21} to the first relation~\eqref{QQPQ} for $\bfq'$ we derive the third relation~\eqref{DefBr3}. The matrix transposition of the second relation~\eqref{QQPQ}, application of
\begin{align} \label{PQtop}
 &P_\bfq^\top=P_{\bfq'}, &&Q_\bfq^\top=Q_{\bfq'}
\end{align}
and~\eqref{PQ21} gives the fifth relation~\eqref{DefBr3}. \qed

\section{Manin matrices of type $C$ and their minors}
\label{secC}

Consider the idempotent% footnote
\footnote{
In the notations of~\cite{Smm} the idempotent~\eqref{Cq} is a deformation of $\wt B_{2r}$. 
}
\begin{align} \label{Cq}
 &C_\bfq=\rho_\bfq(v_1)=\frac{1-P_\bfq}2-\frac{Q_\bfq}{2r}.
\end{align}
Let us write the commutation relations for the quadratic algebras $\Xi_{C_\bfq}(\KK)$ and $\Xi^*_{C_\bfq}(\KK)$. Note that $Q_\bfq C_\bfq=C_\bfq Q_\bfq=0$. By multiplying the matrix relation $(\Psi\otimes\Psi)(1-C_\bfq)=0$ by $Q_\bfq$ from the right we obtain $(\Psi\otimes\Psi)Q_\bfq=0$, so this relation is equivalent to the system
\begin{align} \label{PsiPsiPQ}
 &(\Psi\otimes\Psi)(1+P_\bfq)=0, &&(\Psi\otimes\Psi)Q_\bfq=0.
\end{align}
Analogously, the relations for $\Xi^*_{C_\bfq}(\KK)$ can be written in the form
\begin{align} \label{PQPsiPsi}
 &(1+P_\bfq)(\Psi^*\otimes\Psi^*)=0, &&Q_\bfq(\Psi^*\otimes\Psi^*)=0.
\end{align}
Due to the identities~\eqref{PQtop} the matrix transposition of the relations~\eqref{PQPsiPsi} gives the relations~\eqref{PsiPsiPQ} for $\bfq'=(q_{ij}^{-1},q_i^{-1})$, so the formula $\psi^i\mapsto\psi_i$ gives the isomorphism of quadratic algebras $\Xi^*_{C_\bfq}(\KK)\cong\Xi_{C_{\bfq'}}(\KK)$.

By writing the commutation relations~\eqref{PsiPsiPQ} explicitly we see that $\Xi_{C_\bfq}(\KK)$ is a quadratic algebra with the generators $\psi_{-r},\ldots,\psi_{-1},\psi_1,\ldots,\psi_r$ and relations
\begin{align}
 &\psi_i\psi_j=-q_{ij}\psi_j\psi_i, &&\sum_{i=1}^rq_i\psi_i\psi_{-i}=0.
\end{align}
The first relations define the $\wh q$-Grassmann algebra $\Xi_{A_{\wh q}}(\KK)$ (see~\cite{Manin89}, \cite[\S~3.3]{Smm}), so $\Xi_{C_\bfq}(\KK)$ is the quotient of $\Xi_{A_{\wh q}}(\KK)$ by one additional commutation relation. The quadratic algebra $\Xi^*_{C_\bfq}(\KK)$ is generated by $\psi^{-r},\ldots,\psi^{-1},\psi^1,\ldots,\psi^r$ with the relations
\begin{align}
 &\psi^i\psi^j=-q_{ij}^{-1}\psi^j\psi^i, &&\sum_{i=1}^rq_i^{-1}\psi^i\psi^{-i}=0.
\end{align}

By substituting $A=C_\bfq$ to the first of~\eqref{AXX} we see that $\gX_{C_\bfq}(\KK)$ is the algebra generated by $x^{-r},\ldots,x^{-1},x^1,\ldots,x^r,\lambda$ with the commutation relations $x^ix^j-q_{ji}^{-1}x^jx^i=\varepsilon_iq_i\delta^{i,-j}\lambda$, where $\lambda=\frac1r\sum_{i=1}^r(q_i^{-1}x^ix^{-i}-q_ix^{-i}x^i)$. Again we have the isomorphism $\gX^*_{C_\bfq}(\KK)\cong\gX_{C_{\bfq'}}(\KK)$.

The Molev's definition of {\it Manin matrix of type $C$} (see~\cite[\S~5.6]{MolevSO}, \cite[\S~7.1]{Smm}) is generalised to the case of multi-parametric deformation. We propose the term {\it $(\bfq,\wh p)$-Manin matrix of type $C$} for a $(C_\bfq,A_{\wh p})$-Manin matrix. Let us discuss minors of these matrices.

In~\cite[\S~6.1]{Smm} we calculated $S$-operators for the idempotent $A_{\wh p}$ and $S$-minors for an $(A,A_{\wh p})$-Manin matrix, where $A\in\End(V\otimes V)$ is an arbitrary idempotent (in particular, we can take $A=C_\bfq$).

The $A$-minors of a $(\bfq,\wh p)$-Manin matrix of type $C$ (or, more generally, of a $(C_\bfq,\wt A)$-Manin matrix, where $\wt A\in\End(\wt V\otimes\wt V)$ is an arbitrary idempotent) are expressed through the pairing $A$-operators for the idempotent $C_\bfq$ by the formula~\eqref{MinA}. To find these $A$-operators we remind the following result~\cite[Th.~5.26]{Smm}.

\begin{Lem} \label{ThUkrhosk}
 Consider an idempotent $A\in\End(V\otimes V)$. Let $\rho\colon\gU_k\to\End\big(V^{\otimes k}\big)$ and $\varepsilon\colon\gU_k\to\KK$ be representation and augmentation of an algebra $\gU_k$. Let $s_{(k)}\in\gU_k$ be such that
\begin{align}
 &us_{(k)}=s_{(k)}u=\varepsilon(u)s_{(k)}\qquad\forall\,u\in\gU_k, &&&&\varepsilon(s_{(k)})=1. \label{usk}
\end{align}
Suppose there exist elements $u_1,\ldots,u_{k-1}\in\gU_k$ such that $\rho(u_a)=1-A^{(a,a+1)}$ and $\varepsilon(u_a)=0$ for all $a=1,\ldots,k-1$. If $\rho$ and $\varepsilon$ satisfy
\begin{align}
 (\Psi\otimes\cdots\otimes\Psi)\rho(u) &=\varepsilon(u)(\Psi\otimes\cdots\otimes\Psi), \label{Psipi} \\
 \rho(u)(\Psi^*\otimes\cdots\otimes\Psi^*) &=\varepsilon(u)(\Psi^*\otimes\cdots \otimes\Psi^*) \label{piPsi}
\end{align}
for all $u\in\gU_k$, then $A_{(k)}=\rho(s_{(k)})\in\End\big(V^{\otimes k}\big)$ is the $k$-th $A$-operator for the idempotent $A$.
\end{Lem}

The pairing $A$-operators for the idempotent $A=C_\bfq$ are obtained by application of this fact to the Brauer algebra $\gU_k=\B_k(-2r)$, its representation~\eqref{Rep} and the augmentation
\begin{align}
 &\varepsilon\colon\B_k(\omega)\to\KK, &&\varepsilon(\sigma_a)=1,  &&\varepsilon(\epsilon_a)=0.
\end{align}

\begin{Th}
 The operators
\begin{align} \label{Cqk}
 &C_{\bfq,(k)}=\rho_\bfq(s_{(k)})\in\End(V^{\otimes k}), &&k=2,3,\ldots,r+1,
\end{align}
are the pairing $A$-operators for the idempotent $C_\bfq$.
\end{Th}

\noindent{\bf Proof.} First of all we note that the restriction on $k$ in~\eqref{Cqk} for a fixed $r$ corresponds exactly to the condition~\eqref{omegane} for $\omega=-2r$. The formula~\eqref{sigmask} means exactly that the conditions~\eqref{usk} are valid for the generators $u=\sigma_a$ and $u=\epsilon_a$, so they hold for any $u\in\gU_k$ (the normalisation $\varepsilon(s_{(k)})=1$ is obvious). Since the elements~\eqref{va} satisfy $\rho_\bfq(v_a)=C_\bfq^{(a,a+1)}$ and $\varepsilon(v_a)=1$, $a=1,\ldots,k-1$, we have $\rho_\bfq(u_a)=1-C_\bfq^{(a,a+1)}$ and $\varepsilon(u_a)=0$ for $u_a=1-v_a$. The relations~\eqref{PsiPsiPQ}, \eqref{PQPsiPsi} imply that the conditions~\eqref{Psipi}, \eqref{piPsi} (where $\rho=\rho_\bfq$)
are valid for the generators $u=\sigma_a$, $u=\epsilon_a$ and hence for any element $u\in\gU_k$. Thus Lemma~\ref{ThUkrhosk} can be applied. \qed

Let $M$ be a $(C_\bfq,\wt A)$-Manin matrix, where $\wt A\in\End(\wt V\otimes\wt V)$ is an arbitrary idempotent. Then its $A$-minors are entries of the operator $C_{\bfq,(k)}M^{(1)}\cdots M^{(k)}\in\gR\otimes\Hom(\wt V^{\otimes k},V^{\otimes k})$, where $k\le r+1$. Below we prove that the $k$-th $A$-minors of $M$ vanish for $k\ge r+1$.

\section{Dimensions of components of the quadratic algebra}
\label{secDim}

Here we calculate the dimensions of the graded components $\Xi_{C_\bfq}(\KK)_k$ of the quadratic algebra $\Xi_{C_\bfq}(\KK)$.
For the non-deformed case $q_{ij}=q_i=1$ we found these dimensions in~\cite[\S~7.2]{Smm} by using the abstract (full) trace map $\tr_{1,\ldots,k}\colon\B_k(\omega)\to\KK$ and the formulae~\cite[eq.~(1.52)]{MolevSO}
\begin{align} \label{trsk}
 \tr_{1,\ldots,k}s_{(k)}=\frac{\omega+2k-2}{\omega+k-2}\;{\omega+k-2\choose k}.
\end{align}

\begin{Th}
The graded components of the quadratic algebra $\Xi_{C_\bfq}(\KK)$ have the dimensions
\begin{align}
 &\dim\Xi_{C_\bfq}(\KK)_k=\frac{2r-2k+2}{k}\;{2r+1\choose k-1}, &&k=1,2,\ldots, r+1; \label{dim1} \\
 &\dim\Xi_{C_\bfq}(\KK)_k=0,  &&k\ge r+1. \label{dim2}
\end{align}
\end{Th}

\noindent{\bf Proof.} Recall that the dimension of the vector space $\Xi_{C_\bfq}(\KK)_k$ coincides with the rank of the $k$-th $A$-operator $C_{\bfq,(k)}$ (if it exists)~\cite[eq.~(5.31)]{Smm}. Since the pairing operator is an idempotent, its rank is equal to its trace: $\dim\Xi_{C_\bfq}(\KK)_k=\rk C_{\bfq,(k)}=\tr_{V^{\otimes k}}C_{\bfq,(k)}$. Remind the calculation for the non-deformed case. Due to the commutative diagram~\cite[eq.~(1.77)]{MolevSO} we have $\tr_{V^{\otimes k}}\big(\rho(u)\big)=(-1)^k\tr_{1,\ldots,k}u\;\forall\,u\in\B_k(-2r)$, where $\rho=\rho_\bfq$ for $q_{ij}=q_i=1$. Hence by substituting $\omega=-2r$ to~\eqref{trsk} we obtain the formula~\eqref{dim1} for this case. To prove it for general $\bfq$ we show that $\tr_{V^{\otimes k}}C_{\bfq,(k)}$ does not depend on $\bfq$, if $k=2,\ldots,r+1$ (for $k=1$ it is obvious: $\dim\Xi_{C_\bfq}(\KK)_1=\dim V=2r$). The idempotent~\eqref{sk} can be written in the following form~\cite[Prop.~1.2.5]{MolevSO}:
\begin{align}
 s_{(k)}=h_{(k)}\sum_{t=0}^{[k/2]}\frac{(-1)^t}{2^t t!}{\omega/2+k-2\choose t}^{-1}\sum_{1\le a_i<b_i\le k}\epsilon_{a_1b_1}\cdots\epsilon_{a_tb_t},
\end{align}
where the internal sum is over the non-intersecting pairs of indices $(a_1,b_1),\ldots,(a_t,b_t)$. This expression is well-defined since $\omega/2+k-2=-r+k-2<0$. By using the formulae $\sigma h_{(k)}=h_{(k)}\sigma=h_{(k)}\;\forall\,\sigma\in\SSS_k$ and $\sigma\epsilon_{ab}\sigma^{-1}=\epsilon_{\sigma(a),\sigma(b)}$ we see that the trace $\tr_{V^{\otimes k}}C_{\bfq,(k)}=\tr_{V^{\otimes k}}\rho_\bfq(s_{(k)})$ reduces to a linear combination of the terms $\tr_{V^{\otimes k}}\rho_\bfq(h_{(k)}\epsilon_{12}\cdots\epsilon_{2t-1,2t})$ with coefficients independent of $\bfq$. Let us check that these terms does not depend on $\bfq$ as well.

Note that $\rho_\bfq(h_{(k)})$ equals to the operator $A_{\wh q,(k)}$ with the entries~\eqref{Awhq}, hence we have
\begin{align} \label{trheps}
 \tr_{V^{\otimes k}}\rho_\bfq(h_{(k)}\epsilon_{12}\cdots\epsilon_{2t-1,2t})=\sum_{i_1,\ldots,i_k\atop j_1,\ldots,j_{2t}}(A_{\wh q})^{j_1,\ldots,j_{2t},i_{2t+1},\ldots,i_k}_{i_1,\ldots,i_{2t},i_{2t+1},\ldots,i_k}\prod_{p=1}^t(Q_{\vec q})^{i_{2p-1}i_{2p}}_{j_{2p-1}j_{2p}}.
\end{align}
A term in this sum does not vanish only if the indices $i_1,\ldots,i_k$ are pairwise different, the indices $j_1,\ldots,j_{2t}$ have the form $j_1=i_{\sigma(1)}$,\ldots, $j_{2t}=i_{\sigma(2t)}$ for some $\sigma\in\SSS_{2t}\subset\SSS_k$ and
\begin{align} \label{i2p}
 &i_{2p-1}=-i_{2p}, &&i_{\sigma(2p-1)}=-i_{\sigma(2p)}, &&p=1,\ldots,t.
\end{align}
The conditions~\eqref{i2p} imply that $\sigma=\sigma'\tau$, where $\sigma'=\sigma_{a_1}\cdots\sigma_{a_m}$ for some odd $a_1,\ldots,a_m$ such that $1\le a_1<\ldots<a_m<2t$ and $\tau\in\SSS_{2t}$ permutes the pairs $(1,2), (3,4), \ldots, (2t-1,2t)$ in some way. Hence up to a $\bfq$-independent factor the term in the sum~\eqref{trheps} equals
\begin{align} \label{muprod}
 \mu_I(\wh q,\sigma)\prod_{p=1}^t(q_{i_{2p-1}}q^{-1}_{i_{\sigma(2p-1)}})=
 \mu_I(\wh q,\sigma)\prod_{p=1}^t(q_{i_{2p-1}}q^{-1}_{i_{\sigma'(2p-1)}})=
 \mu_I(\wh q,\sigma)\prod_{s=1}^m q_{i_{a_s}}^2.
\end{align}
Let $\psi_{-r},\ldots,\psi_{-1},\psi_1,\ldots,\psi_r$ be the generators of the $\wh q$-Grassmann algebra $\Xi_{A_{\wh q}}(\KK)$. Due to~\eqref{vecq} they satisfy $\psi_i\psi_{-i}\psi_j\psi_{-j}=\psi_j\psi_{-j}\psi_i\psi_{-i}$, so by taking into account~\eqref{i2p} we obtain
\begin{align}
 \psi_{i_{\sigma(1)}}\cdots\psi_{i_{\sigma(k)}}=
 \psi_{i_{\sigma'\tau(1)}}\cdots\psi_{i_{\sigma'\tau(k)}}=
 \psi_{i_{\sigma'(1)}}\cdots\psi_{i_{\sigma'(k)}}.
\end{align}
By virtue of~\eqref{Lempsisigma2k} this implies $\mu_I(\wh q,\sigma)=\mu_I(\wh q,\sigma')$. The direct calculation of $\mu_I(\wh q,\sigma')=\mu_I(\wh q,\sigma_{a_1}\cdots\sigma_{a_m})$ by the formula~\eqref{muI} gives the product $\prod_{s=1}^m q_{i_{a_s},i_{a_s+1}}$. Hence by using~\eqref{i2p} and \eqref{vecq} again we derive
\begin{align} \label{musimga}
 \mu_I(\wh q,\sigma)=\mu_I(\wh q,\sigma')=\prod_{s=1}^m q_{i_{a_s},i_{a_s+1}}=\prod_{s=1}^m q_{i_{a_s},-i_{a_s}}=\prod_{s=1}^m q_{i_{a_s}}^{-2}.
\end{align}
By substituting~\eqref{musimga} into~\eqref{muprod} we see that the traces~\eqref{trheps} do not depend on $\bfq$. This proves the formula~\eqref{dim1} for general $\bfq$. For $k=r+1$ it gives $\dim\Xi_{C_\bfq}(\KK)_{r+1}=0$ and hence all the higher components $\Xi_{C_\bfq}(\KK)_k$ are also zero. \qed

In particular, the formula~\eqref{dim2} implies that the $k$-th pairing $A$-operator and the corresponding minors vanish for $k\ge r+1$.


\begin{thebibliography}{99}

\bibitem[BW]{BW}
 Birman, J.\,S.; Wenzl, H.: {\it Braids, link polynomials and a new algebra}, Trans. Amer. Math. Soc. \textbf{313} (1989), 249--273.

\bibitem[Br]{Br}
 Brauer, R.: {\it On Algebras Which are Connected with the Semisimple Continuous Groups}, Ann. Math. \textbf{38} (1937), 857--872.

\bibitem[CF]{CF}
 Chervov, A.; Falqui, G.: {\it Manin matrices and Talalaev's formula},
J. Phys. A: Math. Theor. \textbf{41} (2008), No. 19, 194006 (28pp.);
\href{https://arxiv.org/abs/0711.2236}{\tt arXiv:0711.2236}.

\bibitem[CFR]{CFR}
 Chervov, A.; Falqui, G.; Rubtsov, V.: {\it Algebraic properties of Manin matrices~I},  Adv. in Appl. Math.  \textbf{43}  (2009),  no. 3, 239--315;
\eprint{0901.0235}.

\bibitem[CFRS]{CFRS}
 Chervov, A.; Falqui, G.; Rubtsov, V.; Silantyev, A.: {\it Algebraic properties of Manin matrices II: $q$-analogues and integrable systems}, Adv. in Appl. Math. \textbf{60} (2014), 25--89;
\eprint{1210.3529}.

\bibitem[CM]{CM}
 Chervov, A.; Molev, A.: {\it On higher order Sugawara operators},
 IMRN \textbf{2009}, no.9, (2009), 1612--1635; \href{http://arxiv.org/abs/0808.1947}{\tt arXiv:0808.1947}.

\bibitem[FH07]{FH07}
 Foata, D.; Han, G.-N.: {\it A New Proof of the Garoufalidis--L\^{e}--Zeilberger Quantum MacMahon Master Theorem}, J. Algebra \textbf{307} (2007), 424--431; \href{http://arxiv.org/abs/math/0603464}{\tt arXiv:math.CO/0603464}.

\bibitem[FH08]{FH08}
 Foata, D.; Han, G.-N.: {\it A basis for the right quantum algebra and the ``1=q'' principle}, J. Algebraic Combin. \textbf{27} (2008), 163--172; \href{http://arxiv.org/abs/math/0603463}{\tt arXiv:math.CO/0603463}.

\bibitem[GLZ]{GLZ}
 Garoufalidis S.; Le, Thang TQ.; Zeilberger D.: {\it The quantum MacMahon Master Theorem},
Proc. Natl. Acad. of Sci. (PNAS) \textbf{103} (2006), No. 38, 13928--13931;
 \href{http://arxiv.org/abs/math.QA/0303319}{\tt arXiv:math.QA/0303319}.

\bibitem[HS]{HS}
Heckenberger, I.; Sch\"uler, A.: {\it Symmetrizer and Antisymmetrizer of the Birman--Wenzl--Murakami algebras}, Lett.Mat.Phys. \textbf{50} (1999) 45;
 \href{http://arxiv.org/abs/math.QA/0002170}{\tt arXiv:math.QA/0002170}.

\bibitem[I]{I}
 Isaev, A.\,P.: {\it Quantum groups and Yang--Baxter equations}, preprint MPIM (Bonn), MPI 2004-132 (2004).

\bibitem[IM]{IM}
 Isaev, A.\,P.; Molev, A.\,I.: {\it Fusion procedure for the Brauer algebra}, St. Petersburg Math.  J. \textbf{22} (2011), 437--446; \href{http://arxiv.org/abs/0812.4113}{\tt arXiv:0812.4113}.

\bibitem[IMO]{IMO}
 Isaev, A.\,P.; Molev, A.\,I.; Ogievetsky, O.\,V.: {\it A new fusion procedure for the Brauer algebra and evaluation homomorphisms}, Int. Math. Res. Not. \textbf{2012}, issue~11, (2012), 2571--2606.

\bibitem[IO]{IO}
 Isaev, A.\,P.; Ogievetsky, O.\,V.:
{\it Half-quantum linear algebra,} Nankai Series in Pure, Applied Mathematics and Theoretical Physics Symmetries and Groups in Contemporary Physics (2013), 479--486;
\href{http://arxiv.org/abs/1303.3991}{\tt arXiv:1303.3991}.

\bibitem[IR]{IR}
 Isaev, A.\,P.; Rubakov, V.\,A.: {\it Theory Of Groups And Symmetries. Representations of Groups and Lie Algebras, Applications}, World Scientific (2020).

\bibitem[Man87]{Manin87}
 Manin, Yu.\,I.: {\it Some remarks on Koszul algebras and quantum groups},
Ann. de l'Inst. Fourier \textbf{37}, no. 4, (1987), pp. 191--205.

\bibitem[Man88]{Manin88}
 Manin, Yu.\,I.: {\it Quantum Groups and Non-Commutative Geometry},  University of Montreal,
Centre de Recherches Math\'ematiques, Montreal, QC, (1988), 91 pp.

\bibitem[Man89]{Manin89}
 Manin, Yu.\,I.: {\it Multiparametric quantum deformation of the general linear supergroup}, Comm. Math. Phys. \textbf{123}, (1989), 163--175.

\bibitem[Molev]{MolevSO}
 Molev, A.: {\it Sugawara operators for classical Lie algebras}, Providence, RI: American Mathematical Society, (2018), 321 pp.

\bibitem[MR]{MR}
 Molev, A.; Ragoucy, E.: {\it The MacMahon Master Theorem for right quantum superalgebras and higher Sugawara operators for $\mathfrak{gl}(m|n)$}, Moscow Mathematical Journal \textbf{14} (2014), 83--119; \href{https://arxiv.org/abs/0911.3447}{\tt arXiv:0911.3447}.

\bibitem[N]{Nazarov}
 Nazarov, D.: {\it Young’s orthogonal form for Brauer’s centralizer algebra}, J. Algebra \textbf{182} (1996), 664--693.

\bibitem[S20]{Smm}
 Silantyev, A.: {\it Manin matrices for quadratic algebras}, SIGMA \textbf{17}, 066, (2021), 81 pages; \href{https://arxiv.org/abs/2009.05993}{\tt arXiv:2009.05993}.

\bibitem[S21]{Sqrt}
 Silantyev, A.: {\it Quantum Representation Theory and Manin matrices I: finite-di\-men\-sional case},
\href{https://arxiv.org/abs/2108.00269}{\tt arXiv:2108.00269}.

\end{thebibliography}
\end{document}